\documentclass[proceedings,submission%
]{dmtcs}

\usepackage[latin1]{inputenc}

\newlength{\cellsize}
\title[The ABC's of 
affine Grassmannians and Hall-Littlewood polynomials]{The ABC's of 
affine Grassmannians and Hall-Littlewood polynomials}

\author{Avinash J. Dalal
\thanks{Partially supported by the NSF grants DMS--0652641, 
DMS--0638625, DMS--1001898.}
\and 
Jennifer Morse
\thanks{Partially supported by the NSF grants DMS--0652641, DMS--0638652.}}

\address{Department of Mathematics, Drexel University, Philadelphia, PA 19104, U.S.A.}


\usepackage{amsmath,amssymb,latexsym,graphics,color,graphicx}
\usepackage[all]{xy}
\usepackage{tabmac_avi}
\usepackage{color}
\usepackage{algorithmic}

\newtheorem{theorem}{Theorem}

\newtheorem{cor}[theorem]{Corollary}
\newtheorem{lma}[theorem]{Lemma}

\newtheorem{defn}[theorem]{Definition}
\newtheorem{remark}[theorem]{Remark}
\newtheorem{exmpl}[theorem]{Example}
\newtheorem{conj}[theorem]{Conjecture}

\def \la {\lambda}

\def\la {\lambda}

\def\shape{ {\rm {shape}}}

\begin{document}
\newcommand{\apnd}{A}
\newcommand{\rbar}{\bar{R_1}}
\newcommand{\boldx}{\textbf{x}}
\newcommand{\HRule}{\rule{\linewidth}{0.5mm}}
\newcommand{\emptbx}{\mbox{}}
\newcommand{\filbx}{\rule[-0.45mm]{4mm}{4mm}}

\maketitle
\begin{abstract}
\paragraph{Abstract.}
We give a new description of the Pieri rule for $k$-Schur functions 
using the Bruhat order on the affine type-$A$ Weyl group. In
doing so, we prove a new combinatorial formula for representatives 
of the Schubert classes for the cohomology of affine Grassmannians.
We show how new combinatorics involved in our formulas 
gives the Kostka-Foulkes polynomials and discuss how this can be 
applied to study the transition matrices between 
Hall-Littlewood and $k$-Schur functions.  

Nous pr\'esentons une nouvelle description,
issue de l'ordre de Bruhat du groupe
de Weyl affine de type $A$, de la r\`egle de Pieri pour les
fonctions $k$-Schur.  Ce faisant, nous obtenons une nouvelle formule
combinatoire pour les repr\'esentants des classes de Schubert de
la cohomologie des Grassmannienne affines.  
Nous d\'ecrivons aussi comment
notre approche permet d'obtenir les polyn\^omes de Kostka-Foulkes et comment
elle peut \^etre appliqu\'ee \`a l'\'etude des matrices de transition
entre les polyn\^omes de Hall-Littlewood et les fonctions $k$-Schur.
\end{abstract}

\section{Introduction}
\label{sec:in}

The dual $k$-Schur functions arose in \cite{[LMhecke]} where it was 
shown that their coproduct encodes
structure constants of the Verlinde fusion algebra for the 
$\widehat{\mathfrak{sl}}_n$ Wess--Zumino--Witten models 
\cite{TUY:1989,Verlinde:1988} or equivalently,
the 3-point Gromov-Witten invariants of genus zero
(e.g. \cite{[Wi]}).
It was further proven in
\cite{Lam} that these functions represent cohomology classes of the affine 
Grassmannian.  Dual $k$-Schur functions are defined as the weight 
generating functions of $k$-tableaux; a combinatorial object encoding 
successions of saturated chains in weak order on the affine Weyl 
group $\tilde A^{k+1}$ (see \S~\ref{sec:tableauxIntro}).
Here, we find that there is a natural realization of weak 
saturated chains of length $\ell$ as length $k-\ell$ 
saturated chains in 
the strong (Bruhat) order on $\tilde A^{k+1}$.  This
enables us to prove that the dual $k$-Schur functions 
$\mathfrak S_\lambda^{(k)}$ can be described in terms of
a new combinatorial object called affine Bruhat counter-tableau,
or $ABC$, that encodes successions of
strong chains (see Definition~\ref{def:abc}).
We prove, for a fixed positive integer $k$ and $k+1$-core $\lambda$
that
\begin{equation}
\label{abcgenfun}
\mathfrak{S}_{\lambda}^{(k)} = \sum_{A}x^A\,,
\end{equation}
over the $ABC$'s of inner shape $\lambda$.

Our new interpretation for $\{\mathfrak S_\lambda^{(k)}\}_{\lambda_1\leq k}$
puts these cohomology classes in the same combinatorial vein as 
the homology classes.  That is, representatives for the homology 
classes of the affine Grassmannian are given by the $k$-Schur functions 
$s_\lambda^{(k)}$ of \cite{[LMproofs]} which were defined
in terms of chains in the strong order on $\tilde A^{k+1}$ by \cite{[LLMS]}.  
In fact, our work relies on a reformulation of the Pieri rule
for $s_\lambda^{(k)}$; we show that this rule for computing 
the $k$-Schur expansion of
the product $h_r s_\lambda^{(k)}$ can be described by certain
saturated strong chains.

A strong motivation for introducing these new combinatorial ideas is 
to shed light on open problems in the theory of 
Macdonald polynomials.  The origin of $k$-Schur functions was in the idea 
that there exists a more refined basis than Schur functions
upon which Macdonald polynomials $H_\mu[X;q,t]$ expand positively.
In particular, a family of polynomials $A_\lambda^{(k)}[X;t]$
was introduced by \cite{[LLM]} and it was conjectured that,
for $\mu_1\leq k$,
\begin{equation}
\label{mackschur}
H_\mu[X;q,t] = \sum_{\lambda} K_{\lambda\mu}^{(k)}(q,t)\,
A_\lambda^{(k)}[X;t]
\end{equation}
and $K_{\lambda\mu}^{(k)}(q,t)\in\mathbb N[q,t]$.
When $k\geq |\lambda|$, $A_\lambda^{(k)}[X;t]$ reduces to 
the Schur function $s_\lambda$.  In this case, the conjecture
reduces to the conjecture of \cite{[M2]} that $K_{\lambda\mu}^{(\infty)}(q,t)$ 
is a positive sum of monomials in $q$ and $t$.  The polynomiality 
of $K_{\lambda\mu}^{(\infty)}(q,t)$ was proven independently  by
\cite{KiNo96,Sah96,Kno97,GaTe96,LaVi95b} and 
the positivity was finally settled by \cite{Haiman}.  
However, for generic $k$, this remains 
an open problem even in the Hall-Littlewood case $H_\mu[X;0,t]$.

It is conjectured that $A_\mu^{(k)}[X;1]$ are the $k$-Schur functions 
$s_\mu^{(k)}$ prompting our study of $ABC$ combinatorics in the context 
of Hall-Littlewood polynomials.  
To this end, we introduce a simple spin statistic in \S~\ref{sec:hall}
on $ABC$'s and conjecture that the $K^{(k)}_{\lambda\mu}(0,t)$ can be described as
the spin generating function of $ABC$'s.  To provide evidence for our 
conjecture, we prove that the \cite{LSfoulkes} formula for
Kostka-Foulkes polynomials  in terms of charge is equivalent to
$$
K_{\lambda\mu}(t)=\sum_{A} t^{spin(A)} \,,
$$
over all $ABC$'s of weight $\mu$ and inner shape $\lambda$.

\section{$k$-Schur Functions}
\label{sec:tableauxIntro}

Dual $k$-Schur functions are defined as the weight generating functions 
of $k$-tableaux.   To understand the definition of these tableaux,
let us set some notation.  We identify each partition 
$\lambda=(\lambda_1,\ldots, \lambda_n)$ with its
Ferrers shape (having $\lambda_i$ lattice squares in the $i^{th}$ row, 
from the bottom to top).  
A column-strict tableau is a filling of a shape with 
positive integers that weakly decrease along rows and
strictly increase up columns.  A $p$-core is a partition that 
does not contain any cell with hook-length $p$.  We use $\mathcal{C}^p$ 
to denote the set of all $p$-cores.  
The content of cell $(i,j)$ is $j-i$, and its $p$-residue is $j-i$ mod $p$.  
The 5 residues of the 5-core (6,4,3,1,1,1) are 
\begin{displaymath}
\text{\tiny
\tableau[sbY]{0 \cr 1 \cr 2 \cr 3 & 4 & 0 \cr 4 & 0 & 1 & 2 \cr 0 & 1 & 2 & 3 & 4 & 0 }
\,.}
\end{displaymath}
The $p$-degree of a $p$-core $\lambda$,  $deg^p(\lambda)$,
is the number of cells in $\lambda$ whose hook-length is 
smaller than $p$.  In the example above, the core has a $5$-degree
of 11.

Hereafter we work with a fixed integer $k > 0$ and all 
cores (resp. residues) are $k+1$-cores (resp. $k+1$-residues) and
$deg^{k+1}$ will simply be written as $deg$.

\begin{defn}
Let $\nu \in \mathcal{C}^{k+1}$ and let $\alpha = (\alpha_1,\ldots,\alpha_r)$ 
be a composition of $deg(\nu)$. A  \textit{$k$-tableau} of shape $\nu$ and 
$k$-weight $\alpha$ is a column-strict filling of $\nu$ with integers 
$1,2,\ldots,r$ such that the collection of cells filled with letter $i$ are labeled by exactly $\alpha_i$ distinct $k+1$-residues for $1 \leq i \leq r$.  
\end{defn}
The 6-tableaux of weight $(3,3,3,1)$ and shape $(8,2,2) \in \mathcal{C}^7$ are 
\begin{displaymath}
\tableau[sbY]{3_5 & 3_6 \cr 2_6 & 2_0 \cr 1_0 & 1_1 & 1_2 & 2_3 & 3_4 & 3_5 & 3_6 & 4_0}
\hspace{30pt}
\tableau[sbY]{3_5 & 4_6 \cr 2_6 & 3_0 \cr 1_0 & 1_1 & 1_2 & 2_3 & 2_4 & 3_5 & 3_6 & 3_0}
\,.
\end{displaymath}
Note that when $k$ is larger than the biggest hook-length in $\nu$,
a $k$-tableau $T$ of shape $\nu$ and weight $\alpha$ is a column-strict
tableau of weight $\alpha$ since no two diagonals of $T$ 
can have the same residue.  

The dual $k$-Schur functions are the $k$-tableaux generating functions:
for any $\lambda \in \mathcal{C}^{k+1}$, 
\begin{equation}
\label{dualk}
\mathfrak{S}_{\lambda}^{(k)} = 
\sum_{\substack{T: k\text{-tableau} \\\ \shape(T)=\lambda}}
x^{k\text{-}weight(T)}
\,.
\end{equation}

These arose in the context of the quantum cohomology of Grassmannians
in \cite{[LMhecke]}, their definition was then generalized by 
\cite{[Lam]} to give a family of affine Stanley symmetric functions,
and it was shown by \cite{Lam} that they represent Schubert cohomology
classes of affine Grassmannians.  
The term affine Schur function is also used for dual $k$-Schur function.

\subsection{Weak $k$-Pieri Rule}
\label{subsec:weak k-pieri rule}

Our new formulation for $\mathfrak S^{(k)}_\lambda$ in \eqref{abcgenfun}
is derived along the same lines that led to the introduction of $k$-tableaux 
in \cite{[LMcore]}.  This route starts with a family of symmetric functions 
related to $\mathfrak S^{(k)}_\lambda$ by duality.  To be precise, the set of 
$\{\mathfrak S_\lambda^{(k)}\}_{\lambda\in\mathcal C^{k+1}}$ 
forms a basis for 
\begin{equation}
\Lambda/\mathcal I\qquad\text{where}\quad \mathcal I=\langle m_\lambda
\rangle_{\lambda_1>k}\,.
\end{equation}
This quotient is naturally 
paired with the subring $\Lambda^k = \mathbb{Z}[h_1,\ldots,h_k]$ 
of symmetric functions under the Hall-inner product,
\[
\langle h_\lambda,m_\mu\rangle = \delta_{\lambda\mu}
\,.
\]
The basis dual to the set of 
$\{\mathfrak S_\lambda^{(k)}\}_{\lambda\in\mathcal C^{k+1}}$ 
turns out to be the $k$-Schur function basis
$\{s_\lambda^{(k)}\}_{\lambda\in\mathcal C^{k+1}}$, functions
conjectured to be the $t=1$ case of the atoms 
that arose in the Macdonald polynomial study of \cite{[LLM]}.  
By the duality of $\{h_\lambda\}$ and $\{m_\lambda\}$,
 \eqref{dualk} implies that the close examination of 
the expansion
\begin{equation}
h_\mu = \sum_\lambda K_{\lambda\mu}^{(k)} \, s_\lambda^{(k)}
\end{equation}
should reveal the $k$-tableaux.
In fact, starting with $s_\emptyset^{(k)}=1$,
$k$-tableaux are precisely the
objects that encode the combinatorial Pieri rule
for $k$-Schur functions given in \cite{[LMproofs]}.
That is, for $\lambda\in \mathcal{C}^{k+1}$ and $0<\ell\leq k$, 
\begin{displaymath}
h_{\ell}s_{\lambda}^{(k)} = \sum_{\nu \in {H}_{\lambda,\ell}^{(k)}} s_{\nu}^{(k)}
\,,
\end{displaymath}
over a specified set of cores ${H}_{\lambda,\ell}^{(k)}$
described by certain saturated chains in the weak order 
on the affine Weyl group $\tilde A^{k+1}$.  Recall that weak order 
on $\tilde A^{k+1}$ can be realized on cores by the covering relation
\begin{displaymath}
\beta \lessdot \nu \Longleftrightarrow \nu = \beta + \text{all addable corners of one fixed residue}.
\end{displaymath}
\begin{defn}
For $0<\ell \leq k$ and $k+1$-cores $\lambda$ and $\nu$, the skew shape $\nu/\lambda$
is a \emph{weak $\ell$-strip} if there is a weak saturated chain of cores
\begin{displaymath}
\lambda = \gamma^0 \lessdot \gamma^1 \lessdot \cdots \lessdot \gamma^{\ell} = 
\nu\,,
\end{displaymath} 
where $\nu/\lambda$ is 
a horizontal strip and the rightmost cell 
of $\gamma^i/\gamma^{i-1}$ is to the left of that in $\gamma^{i+1}/\gamma^i$ for all $i = 1, \ldots, \ell-1$.    
\end{defn}
It should be observed that if there exists in the previous definition such a weak saturated chain then it is unique.
\begin{exmpl}
\label{weak2}
The skew shape $(4,1,1,1)/(3,1,1)$ of 4-cores is a weak 2-strip 
as there is the saturated chain
\begin{displaymath}
\tiny\tableau[sbY]{ \cr \cr & & } \lessdot \tableau[sbY]{ \cr \cr \cr & & } \lessdot \tableau[sbY]{ \cr \cr \cr & & & }\,.
\end{displaymath}
\end{exmpl}
The terms appearing in the summand of the $k$-Pieri rule are 
simply the set of weak $\ell$-strips.  That is,
\begin{displaymath}
{H}_{\lambda,\ell}^{(k)} = \{ \nu : \nu/\lambda \text{ is a weak } \ell\text{-strip} \}. 
\end{displaymath}

\subsection{Strong $k$-Pieri Rule}
\label{sec:strong k-pieri rule}
Recall that the
strong (Bruhat) order on affine type $A$ is realized on cores 
by the covering relation:
\begin{displaymath}
  \rho \lessdot_B \gamma \Longleftrightarrow \rho \subseteq \gamma  
\text{ and } deg(\gamma) = deg(\rho)+1.
\end{displaymath}
An important fact about strong covers is useful in our study.
An \textit{$n$-ribbon} $R$ is a skew diagram $\la/\mu$ 
consisting of $n$ rookwise connected cells such that there is
no $2\times 2$ shape contained in $R$.
We refer to the southeasternmost cell of a ribbon as its
\textit{head}.

\begin{lma}\cite{[LLMS]}
Let $\rho\lessdot_B\gamma$ be cores.  Then
\begin{enumerate}
\item Each connected component of $\rho/\gamma$ is a ribbon. 
\item The components are translates of each other and their heads
have the same residue. 
\end{enumerate}
\end{lma}

\begin{defn}
  $(\gamma,c)$ is a marked strong cover of a $k+1$-core $\rho$ if $\rho \lessdot_B \gamma$ and $c$ is the content of the head of a ribbon in $\gamma/\rho$.
\end{defn}

\begin{exmpl}
Consider the 4-cores $\rho=(4,1,1,1),\gamma=(6,3,1,1)$.  
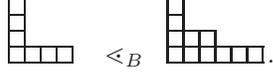
\begin{displaymath}
\tiny\tableau[sbY]{ \cr \cr \cr & & & } \hspace{10pt} \lessdot_B \hspace{7pt} \tableau[sbY]{ \cr \cr & & \cr & & & & & }\,.
\end{displaymath}
$(\gamma,5)$ is a marked strong cover of $\rho$, as is $(\gamma,1)$.
\end{exmpl}

\begin{defn}
For $0 < \ell \leq k$ and $k+1$-cores $\lambda$ and $\gamma$, a strong $\ell$-strip from $\lambda$ to $\gamma$ is a saturated chain in the strong order
\begin{displaymath}
\lambda = \gamma^0 \lessdot_B \gamma^1 \lessdot_B \ldots \lessdot_B \gamma^{\ell} = \gamma
\end{displaymath}
together with a content vector $c = (c_1,c_2,\cdots,c_{\ell})$ such that
\begin{enumerate}
\item $(\gamma^i,c_i)$ is a marked strong cover of $\gamma^{i-1}$,
\item $c_{i-1} < c_i$ for $2 \leq i \leq \ell$.
\end{enumerate}
\end{defn}

\begin{exmpl}
The following sequence of marked strong covers is a strong 2-strip from $(1)$ to $(3,1)$
\begin{displaymath}
\tiny\tableau[sbY]{ & \bl} \hspace{5pt} \lessdot_B \hspace{5pt} \tiny\tableau[sbY]{ & \bullet} \hspace{5pt} \lessdot_B \hspace{5pt} 
\tiny\tableau[sbY]{ \cr & & \bullet}
\end{displaymath} 
since the content vector $c = (1,2)$ is increasing.
\end{exmpl}

The strong $\ell$-strips help to encode the combinatorial strong 
Pieri rule for the dual $k$-Schur functions.  That is, for 
$\lambda \in \mathcal{C}^{k+1}$ and $0 < \ell \leq k$, 
\begin{displaymath}
h_{\ell} \mathfrak{S}_{\lambda}^{(k)} = \sum_{\gamma} d_{\gamma} \mathfrak{S}_{\gamma}^{(k)},
\end{displaymath}
where $d_{\gamma}$ is the number of strong $\ell$-strips from $\lambda$ to $\gamma$.

\section{Bottom Strong Strips}
\label{sec:kPieri}

We have discovered that there is a natural association between
weak strips and certain strong strips.
In this association, a saturated chain in weak order of
length $\ell$ is identified with a saturated chain in strong order
of length $k-\ell$.

Our notion comes from a close examination 
of the core $(k+\lambda_1,\lambda)$ constructed from any
$\lambda \in \mathcal{C}^{k+1}$.

\begin{defn}
For $0 < \ell \leq k$ and $k+1$-cores $\lambda$ and $\nu$, the skew shape $(k+\lambda_1,\lambda)/\nu$ 
is a \emph{bottom strong $\ell$-strip} if there is a saturated chain 
of cores
\begin{displaymath}
\nu = \nu^0 \lessdot_B \nu^1 \lessdot_B \cdots 
\lessdot_B \nu^{k-\ell}=(k+\lambda_1,\lambda)\,,
\end{displaymath}
where 
\begin{enumerate}
 \item $(k+\lambda_1,\lambda)/\nu$ is a horizontal strip.  
 \item The bottom rightmost cell of $\nu^i$ (which is always in the first row) is also a cell in $\nu^i/\nu^{i-1}$, for $1 \leq i \leq k-\ell$.
\end{enumerate}
\end{defn}
It turns out that in the previous definition if such a chain exists then it is unique.
\begin{remark}
Every bottom strong $\ell$-strip $(k+\lambda_1,\lambda)/\nu$ is a strong $\ell$-strip if we take the content of $\nu^i$ to be the content of its bottom rightmost cell in the corresponding saturated chain of cores
\begin{displaymath}
\nu = \nu^0 \lessdot_B \nu^1 \lessdot_B \cdots 
\lessdot_B \nu^{k-\ell}=(k+\lambda_1,\lambda)\,.
\end{displaymath}
\end{remark}

\begin{exmpl}
The skew shape $(8,3)/(4,2)$ of 6-cores is a bottom strong 3-strip as there is the 
saturated chain
\begin{displaymath}
\tiny\tableau[sbY]{&\cr&&&}\lessdot_B
\tableau[sbY]{&\cr&&&&&&}
\lessdot_B
\tableau[sbY]{&&\cr&&&&&&&}.
\end{displaymath}
\end{exmpl}

\begin{exmpl}
The skew shape $(6,3,1,1)/(4,1,1,1)$ of 
4-cores is a bottom strong 2-strip as there is the saturated chain 
\begin{displaymath}
\tiny\tableau[sbY]{ \cr \cr \cr & & & } \hspace{10pt} \lessdot_B \hspace{7pt} \tableau[sbY]{ \cr \cr & & \cr & & & & & }\,.
\end{displaymath}
\end{exmpl}

This is an example of a bottom strong $2$-strip which is
also a weak $2$-strip (compare to Example~\ref{weak2}).
In fact, there is a close connection between bottom $\ell$-strips and weak $\ell$-strips.

\begin{theorem}\label{weakStrongConnection}
For $0<\ell\leq k$ and $k+1$-cores $\lambda$ and $\nu$, 
\begin{displaymath}
  (k+\lambda_1,\lambda)/\nu \text{ is a bottom strong } \ell\text{-strip} \Longleftrightarrow \nu/\lambda \text{ is a weak } \ell\text{-strip}.
\end{displaymath}
\end{theorem}
One immediate application of this result is that the $k$-Pieri rule 
can instead be given in terms of bottom strong $\ell$-strips.
\begin{cor}\label{strongkpieri}
For any $0<\ell\leq k$ and $\lambda\in \mathcal{C}^{k+1}$,
\begin{displaymath}
h_\ell \, s_\lambda^{(k)} = \sum_{\substack{\nu\in\mathcal C^{k+1} \\\
(k+\lambda_1,\lambda)/\nu: \text{ bottom strong }\ell\text{-strip}}}
s_\nu^{(k)}.
\end{displaymath}
\end{cor}

\begin{exmpl}
When $k=3$, the set of $\nu$ such that $(6,3,1,1)/\nu$ is a bottom strong 2-strip is
\begin{displaymath}
\{ (3,3,1,1), (5,2,1), (4,1,1,1) \}\,,
\end{displaymath}
and thus
\begin{displaymath}
  h_{2}s_{(3,1,1)}^{(3)} = s_{(3,3,1,1)}^{(3)} + s_{(5,2,1)}^{(3)} + 
s_{(4,1,1,1)}^{(3)}.
\end{displaymath}
\end{exmpl}

\section{Affine Bruhat Counter-Tableaux}
\label{sec:abc}

The iteration of Corollary~\ref{strongkpieri} gives rise to 
a new combinatorial structure whose enumeration describes
the coefficients in
\begin{equation}
\label{hinkschur}
h_{\lambda} = \sum_{\nu}K_{\nu\lambda}^{(k)} s_{\nu}^{(k)}\,.
\end{equation}
To describe this structure, first let us set more notations.
Recall that a counter-tableau $A$ is a filling of a skew
shape with numbers that strictly decrease up columns and weakly decrease along the rows.
Given a counter-tableau $A$, let $A^{(x)}$ denote the subtableau made up of 
the rows of $A$ weakly higher than row $x$ (where row 1 is this time the \textit{highest} row).
Let $A_{>i}$ denote the restriction of $A$ to letters strictly larger
than $i$ where empty cells in a skew are considered to contain
$\infty$.

\begin{defn}
\label{def:abc}
For a composition $\alpha$ whose entries are not larger than $k$,
a skew counter-tableau $A$ is an \textit{affine Bruhat counter-tableau} 
(or $ABC$) of $k$-weight $\alpha$ if
$$(k+\lambda_1^{(i-1)},\lambda^{(i-1)})/\lambda^{(i)}\;\;\text{is a bottom strong $\alpha_i$-strip
\,for all $1\leq i\leq\ell(\alpha)$}\,,
$$
where
$\lambda^{(x)}=\shape(A^{(x)}_{>x})$.
The inner shape of $A$ is $\lambda^{(\ell(\alpha))}$.
\end{defn}
In the construction of a counter-tableau, we start with the
empty shape $\lambda^{(0)}$, successively adding bottom strong strips
and tiling 
$(k+\lambda^{(i-1)}_1,\lambda^{(i-1)})/\lambda^{(i)}$ with $i$-ribbons
at each step.  
\begin{exmpl}
With $k=5$, an $ABC$ of 5-weight $(3,3,1)$ is
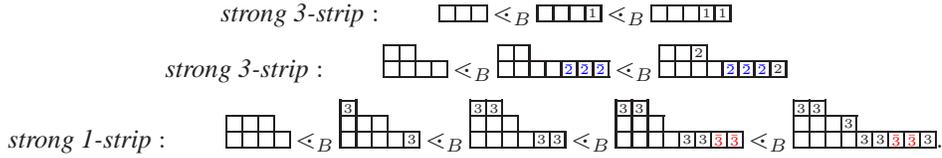
\begin{displaymath}
\tableau[sbY]{3&3&2&1&1\cr&&&3&{\color{blue}\bar 2}&{\color{blue}\bar 2}&{\color{blue}\bar 2}&2 \cr &&&&3&3&{\color{red}\bar 3}&{\color{red}\bar 3}&3},
\end{displaymath}
since
\begin{displaymath}
\text{strong 3-strip}: 
\qquad
\tiny\tableau[sbY]{&&}\lessdot_B \tableau[sbY]{&&&1}
\lessdot_B \tableau[sbY]{&&&1&1}
\end{displaymath}
\begin{displaymath}
\text{strong 3-strip}: 
\qquad
\tiny\tableau[sbY]{&\cr&&&}\lessdot_B
\tableau[sbY]{&\cr&&&&{\color{blue}\bar 2}&{\color{blue}\bar 2}&{\color{blue}\bar 2}}
\lessdot_B
\tableau[sbY]{&&2\cr&&&&{\color{blue}\bar 2}&{\color{blue}\bar 2}&{\color{blue}\bar 2}&2}
\end{displaymath}
\begin{displaymath}
\text{strong 1-strip}: 
\qquad
\tiny\tableau[sbY]{&&\cr&&&}
\lessdot_B
\tiny\tableau[sbY]{3\cr&&\cr &&&&3}
\lessdot_B
\tiny\tableau[sbY]{3&3\cr&&\cr&&&&3&3}
\lessdot_B
\tiny\tableau[sbY]{3&3\cr&&\cr&&&&3&3&\color{red}\bar 3 & \color{red}\bar 3}
\lessdot_B
\tiny\tableau[sbY]{3&3\cr&&&3\cr&&&&3&3&\color{red}\bar 3&\color{red}\bar 3 & 3}.
\end{displaymath}
The black letters are ribbons of size one, red letters make a ribbon of size two and blue letters make a ribbon of size 3 (or for those with out color the ribbons are depicted with a bar).
\end{exmpl}

For a composition $\alpha$ whose entries are not larger than $k$, there is a known bijection between the set of all standard $k$-tableaux of $k$-weight $\alpha$ and the set of all $ABC$'s of $k$-weight $\alpha$.  Letting $K_{\nu\lambda}^{(k)}$ denote the number of affine Bruhat 
counter-tableaux of inner shape $\nu$ and 
$k$-weight $\lambda$, it is not hard to see that it is this
set of $ABC$'s being enumerated by the coefficients in \eqref{hinkschur}
and by duality we can prove that the dual $k$-Schur functions
are the $ABC$ generating functions.
\begin{theorem}
For any $\lambda \in \mathcal{C}^{k+1}$,
\begin{displaymath}
\mathfrak{S}_{\lambda}^{(k)} = \sum_{A} x^{A}
\end{displaymath}
where the sum is over all affine Bruhat counter-tableaux of inner 
shape $\lambda$, and $x^{A} = x^{k\text{-}weight(A)}$.
\end{theorem}

\section{Hall-Littlewood Expansions}
\label{sec:hall}
A motivation for Theorem~\ref{weakStrongConnection} is in its 
application to the study of the integral form of Macdonald 
polynomials (e.g. \cite{Macbook,Ber09}) and the $k$-Schur
expansion coefficients $K_{\lambda\mu}^{(k)}(q,t)$
in \eqref{mackschur}.  For $\mu\vdash n$, it is known that
the Macdonald polynomial $H_\mu[X;q,t]$
reduces to $h_{1^n}$ when $q=t=1$.  
Thus, $$K_{\lambda\mu}^{(k)}(1,1)=K_{\lambda 1^n}^{(k)}$$
are the coefficients in \eqref{hinkschur} when $\mu=(1^n)$. We have
seen that these count the number of $ABC$'s of shape $\lambda$
with $k$-weight $(1^n)$ (called \textit{standard} $ABC$'s).
Assuming the positivity of $K_{\lambda\mu}^{(k)}(q,t)$, this
leads us to believe there exists a pair of statistics 
(non-negative integers)  $a(A)$ and $b(A)$ associated to 
each $ABC$ so that
\begin{displaymath}
K_{\lambda\mu}^{(k)}(q,t) = 
\sum_{\substack{A=\text{standard ABC} \\\ \text{inner}~\shape(A)=\lambda}} 
t^{a_\mu(A)} q^{b_\mu(A)}
\,.
\end{displaymath}

As a special case, when $k=n$, these statistics would give a combinatorial 
formula for the $q,t$-Kostka polynomials arising in the Schur expansion 
of the Macdonald polynomials -- a long-standing open problem in the field.  
Setting $q=0$, the Macdonald polynomials reduce to Hall-Littlewood polynomials.  
In fact, when $k=n$, \cite{LSfoulkes} give a beautiful combinatorial 
formula for the Kostka-Foulkes polynomials $K_{\lambda\mu}^{(n)}(0,t)$ by 
defining a statistic on column-strict tableaux.
It is interesting to note that $k$-tableaux reduce to
column-strict tableaux when $k=n$, whereas the $ABC$'s are
truly of a different spirit even when $k=n$.
We thus start by introducing a statistic on $ABC$'s when
$k=n$.  

Here, we will restrict our attention to \textit{standard}
$ABC$'s -- those of weight $(1^n)$.  It turns out that a standard 
$ABC$ has only 1 or 2-ribbons.  In fact, when $k=n$, there is a simple 
statistic that gives the Kostka-Foulkes polynomials.

\begin{defn}
\label{stdSpinkinf}
For each $ABC$ $A$ of n-weight $1^n$, let
$$
spin(A)=\sum_i i\,\chi(\text{row $i$ has a 2-ribbon and
it is not east of a 2-ribbon in row $i+1$})\,,
$$
where row 1 is the bottommost row of $A$, and 
where $\chi$ is the indicator function evaluating to
zero when the statement is false and to 1 otherwise.  Recall that a ribbon is east of another ribbon if the head of the former is to the east of the latter's.
\end{defn}

\begin{exmpl}
Consider the $ABC$ of 5-weight $1^5$
\begin{displaymath}
         A = \tableau[scY]{2&1&1&1&1 \cr 4&3&2&2&2&2 \cr
&5&3&3&{\color{red}\bar 3}&{\color{red}\bar 3}&3 \cr
&&4&4&4&{\color{red}\bar 4}&{\color{red}\bar 4} \cr
&&5&{\color{red}\bar 5}&{\color{red}\bar 5}&5&5}.
\end{displaymath}
2-ribbons are colored (or for those without color, 2-ribbons 
are depicted with a bar).
Of these 3 2-ribbons, only the one in row 1 and row 3 contribute to the spin
since the 2-ribbon in row 2 is east of that in row 3.  The
spin is thus 1+3=4.
\end{exmpl}

\begin{theorem}\label{stdABCcch}
For $\lambda\vdash n$,
\begin{equation}
K_{\lambda, 1^n}(0,t)=
K_{\lambda 1^n}^{(n)}(0,t)= \sum_{\substack{A: \text{ $ABC$ of $n$-weight $1^n$} \\\
\text{inner shape}(A)=\lambda}} t^{spin(A)}
\,.
\end{equation}
\end{theorem}

\begin{exmpl}
The set of all $ABC$'s of 3-weight $1^3$ are
\begin{displaymath}
 \left\{
 \hspace{5pt}
 \tableau[scY]{ & 1 & 1 \cr & 2 &{\color{red}\bar 2}&{\color{red}\bar 2} \cr &
{\color{red}\bar 3} & {\color{red}\bar 3} & 3 }
 \hspace{10pt}
 \tableau[scY]{3&1&1 \cr & 2 &{\color{red}\bar 2}&{\color{red}\bar 2} \cr &&3&3 }
 \hspace{10pt}
 \tableau[scY]{2&1&1 \cr & 3 & 2 & 2 \cr &&{\color{red}\bar 3}&{\color{red}\bar 3}&3 }
 \hspace{10pt}
 \tableau[scY]{2&1&1 \cr 3 & 3 & 2 & 2 \cr & & & 3 & 3 }
 \hspace{5pt}
 \right\}
\end{displaymath}
The spins of these $ABC$'s, respectively, are \{3, 2, 1, 0\}.
We thus have that
\begin{align*}
 H_{1^3}[X,t] &= t^3\,s_{(1,1,1)}+ (t^2+t)\,s_{(2,1)} +
s_{(3)}\,.
\end{align*}
\end{exmpl}

The beauty of this formulation is that it gives a hint
into finding a statistic to describe the more general Hall-Littlewood 
coefficients $K_{\lambda,1^n}^{(k)}(t)$.
To do this, we need only the extra concept of an
\textit{offset} in an $ABC$. 
In an $ABC$, for $i>1$, an $i$-ribbon $R$ 
is an \textit{offset} if there is a lower $i$-ribbon filled with 
the same letter as $R$ whose head has the same residue 
as the head of $R$.  The number of offsets in the $ABC$ is
denoted \textit{off}$(ABC)$.
\begin{exmpl}
Consider the $ABC$ of $3$-weight $1^5$.
\begin{displaymath}
       A = \tableau[scY]{3&1&1 \cr  & 2 & {\color{red}\bar 2} & {\color{red}\bar 2}
\cr &4&3&3 \cr  &{\color{red}\bar 5}&{\color{red}\bar 5}&4&4 \cr & & & 5 &
{\color{red}\bar 5} & {\color{red}\bar 5} }
\end{displaymath}
Here we see that $A$ has only one offset
$\tableau[scY]{{\color{red}\bar 5}&{\color{red}\bar 5}}$ in the second row from
the bottom.  Observe that it is an offset since there is another 2-ribbon of the same residue in the first row with the same letter.  This tells us that off$(A) = 1$.
\end{exmpl}

\begin{defn}
Let $A$ be an $ABC$ of $k$-weight $1^n$ and define
\begin{displaymath}
spin^k(A) = \text{off}(A) + \sum_{i} i\,\chi(\text{row $i$ has a 2-ribbon 
and it is not east of a 2-ribbon in row $i+1$})\,.
\end{displaymath}
\end{defn}

When $k=n$, an $ABC$ never has any offsets and thus the above
definition reduces to Definition~\ref{stdSpinkinf}.

\begin{conj}\label{stdtKostkaSpinAbcAtk}
For any $k+1$-core $\lambda$ where $deg(\lambda)=n$,
\begin{equation}
K_{\lambda 1^n}^{(k)}(0,t)= \sum_{\substack{A: \text{ $ABC$ of weight $1^n$} \\\
\text{inner shape}(A)=\lambda}} t^{spin^k(A)}
\,.
\end{equation}
\end{conj}

\begin{exmpl}
The set of all $ABC's$ of 2-weight $1^4$ are
\begin{displaymath}
 \left\{
 \hspace{5pt}
 \tableau[scY]{ & 1 \cr & {\color{red}\bar 2}&{\color{red}\bar 2} \cr & & 3
\cr & & {\color{red}\bar 4}&{\color{red}\bar 4} }
 \hspace{10pt}
 \tableau[scY]{4& 1 \cr & {\color{red}\bar 2}&{\color{red}\bar 2} \cr &4& 3 \cr & & &4 }
 \hspace{10pt}
 \tableau[scY]{2& 1 \cr & 3 & 2 \cr &{\color{red}\bar 4}&{\color{red}\bar 4}&3
\cr & & &{\color{red}\bar 4}&{\color{red}\bar 4} }
 \hspace{10pt}
 \tableau[scY]{2& 1 \cr4& 3 & 2 \cr & & 4 & 3 \cr & & & & 4 }
 \hspace{5pt}
 \right\}
\end{displaymath}
Hence, our conjecture checks out against the expansion of $H_{1^4}[X,t]$ 
in terms of 2-Schur functions:
\begin{align*}
H_{1^4}[X,t] &= s_{(2,2,1,1)}^{(2)} + t^2s_{(3,1,1)}^{(2)} +
t^3s_{(3,1,1)}^{(2)} + t^4s_{(4,2)}^{(2)}\,.
\end{align*}
\end{exmpl}

\section{Related and Future work}

Central to the proof of the \cite{LSfoulkes} formula
for Hall-Littlewood polynomials is the {\it plactic monoid}, \cite{[LS]}.
The {plactic monoid} is the quotient of the free 
monoid $A^*$ on the totally ordered set $A$ by 
{\it Knuth relations}, defined on words so $w\sim w'$ 
if and only if they are sent to the same 
tableau under RSK-insertion (see \cite{Knu70}).  This establishes a 
bijection $T \mapsto [w] \in (A^* / \sim)$ and the monoid
can be viewed as the set of tableaux on letters in 
$A$ with multiplication defined by insertion.  
The ring of symmetric functions $\Lambda$ can be identified with 
a subring of the monoid algebra $\mathbb{Z}[q,t][A^* /\sim]$
by sending the Schur function $s_\lambda$ to the sum of all tableaux 
with shape $\lambda$.  Although the plactic monoid is noncommutative,
the combinatorial nature of computation is more evident.  

A refinement of the plactic monoid to a structure on $k$-tableaux 
that can be applied to combinatorial problems involving $k$-Schur 
functions is partially given in \cite{[LLMS2]} by a bijection
compatible with the RSK-bijection.  A deeper understanding of this
highly intricate bijection is underway.  Towards this effort,
Lapointe and Pinto have discovered a statistic on $k$-tableaux
that is compatible with the bijection \cite{[LLMEP]}.
We conjecture that their statistic matches the spin of our $ABC$'s 
and are working to put $ABC$'s in a context that simplifies 
the bijection.

\bibliographystyle{alpha}
\bibliography{ktheoryjune}

\begin{thebibliography}{LLMS10}

\bibitem[Ber09]{Ber09}
Fran{\c{c}}ois Bergeron.
\newblock {\em Algebraic combinatorics and coinvariant spaces}.
\newblock CMS Treatises in Mathematics. Canadian Mathematical Society, Ottawa,
  ON, 2009.

\bibitem[GT96]{GaTe96}
A.~M. Garsia and G.~Tesler.
\newblock Plethystic formulas for {M}acdonald $q,t$-{K}ostka coefficients.
\newblock {\em Adv. Math.}, 123:144--222, 1996.

\bibitem[Hai01]{Haiman}
M.~Haiman.
\newblock Hilbert schemes, polygraphs, and the macdonald positivity conjecture.
\newblock {\em J. Amer. Math. Soc.}, 14(4):941--1006, 2001.

\bibitem[KN98]{KiNo96}
A.~N. Kirillov and M.~Noumi.
\newblock Affine hecke algebras and raising operators for macdonald
  polynomials.
\newblock {\em Duke Math. J.}, 93(1):1--39, 1998.

\bibitem[Kno97]{Kno97}
F.~Knop.
\newblock Integrality of two variable {K}ostka functions.
\newblock {\em J. Reine Angew. Math.}, 482:177--189, 1997.

\bibitem[Knu70]{Knu70}
D.~E. Knuth.
\newblock Permutations, matrices, and generalized {Y}oung tableaux.
\newblock {\em Pacific J. Math.}, 34:709--727, 1970.

\bibitem[Lam06]{[Lam]}
T.~Lam.
\newblock Affine stanley symmetric functions.
\newblock {\em Amer. J. Math.}, 128(6):1553--1586, 2006.

\bibitem[Lam08]{Lam}
T.~Lam.
\newblock Schubert polynomials for the affine grassmannian.
\newblock {\em J. of Amer. Math. Soc.}, 21(1):259--281, 2008.

\bibitem[LLM03]{[LLM]}
L.~Lapointe, A.~Lascoux, and J.~Morse.
\newblock Tableau atoms and a new macdonald positivity conjecture.
\newblock {\em Duke Math. J.}, 116(1):103--146, 2003.

\bibitem[LLMS]{[LLMS2]}
T.~Lam, L.~Lapointe, J.~Morse, and M.~Shimozono.
\newblock The poset of $k$-shape and branching of $k$-schur functions.
\newblock {\em to be published in Mem. Amer. Math. Soc.}

\bibitem[LLMS10]{[LLMS]}
T.~Lam, L.~Lapointe, J.~Morse, and M.~Shimozono.
\newblock Affine insertion and pieri rules for the affine grassmannian.
\newblock {\em Mem. Amer. Math Soc.}, 208(977), 2010.

\bibitem[LM05]{[LMcore]}
L.~Lapointe and J.~Morse.
\newblock Tableaux on $k+1$-cores, reduced words for affine permutations, and
  $k$-schur function expansions.
\newblock {\em J. Combin. Theory Ser.}, 112:44--81, 2005.

\bibitem[LM07]{[LMproofs]}
L.~Lapointe and J.~Morse.
\newblock A $k$-tableaux characterization for $k$-schur functions.
\newblock {\em Adv. Math.}, pages 183--204, 2007.

\bibitem[LM08]{[LMhecke]}
L.~Lapointe and J.~Morse.
\newblock Quantum cohomology and the $k$-schur basis.
\newblock {\em Trans. Amer. Math. Soc.}, 360:2021--2040, 2008.

\bibitem[LP]{[LLMEP]}
L.~Lapointe and M.~E. Pinto.
\newblock Private communication.

\bibitem[LS78]{LSfoulkes}
A.~Lascoux and M.-P. Sch\"utzenberger.
\newblock Sur une conjecture de h.o. foulkes.
\newblock {\em C.R. Acad. Sc. Paris}, 294:323--324, 1978.

\bibitem[LS81]{[LS]}
A.~Lascoux and M.-P. Sch\"utzenberger.
\newblock Le mono\"\i de plaxique.
\newblock {\em Quaderni della Ricerca scientifica}, 109:129--156, 1981.

\bibitem[LV97]{LaVi95b}
L.~Lapointe and L.~Vinet.
\newblock A rodrigues formula for the macdonald polynomials.
\newblock {\em Adv. Math.}, 130:261--279, 1997.

\bibitem[Mac88]{[M2]}
I.~G. Macdonald.
\newblock A new class of symmetric functions.
\newblock {\em S\'eminaire Lotharingien de Combinatoire}, {\bf B20a} 41pp,
  1988.

\bibitem[Mac95]{Macbook}
I.~G. Macdonald.
\newblock {\em Symmetric functions and Hall polynomials}.
\newblock Clarendon Press, Oxford, 2nd edition, 1995.

\bibitem[Sah96]{Sah96}
S.~Sahi.
\newblock Interpolation, integrality, and a generalization of {M}acdonald's
  polynomials.
\newblock {\em Internat. Math. Res. Notices}, pages 457--471, 1996.

\bibitem[TUY89]{TUY:1989}
A.~Tsuchiya, K.~Ueno, and Y.~Yamada.
\newblock Conformal field theory on universal family of stable curves with
  gauge symmetries.
\newblock {\em Adv. Stud. Pure Math.}, 19:459--566, 1989.

\bibitem[Ver88]{Verlinde:1988}
E.~Verlinde.
\newblock Fusion rules and modular transformations in 2d conformal field
  theory.
\newblock {\em Nuclear Phys. B}, 300(3):360--376, 1988.

\bibitem[Wit95]{[Wi]}
E.~Witten.
\newblock The verlinde algebra and the cohomology of the grassmanian,
  ``geometry, topology and physics".
\newblock {\em Conf. Proc. Lecture Notes Geom. Topology, IV}, pages 357--422,
  1995.

\end{thebibliography}
\label{sec:biblio}

\end{document}